\documentclass[12pt]{article}
\usepackage{epsfig}
\usepackage{amsmath,amssymb}
\newcommand{\R}{I \! \! R}

\newcommand{\C}{I \! \! \! \! {C}}

\newcommand{\uu}{{\underline u}}

\newcommand{\us}{{\underline s}}
\newcommand{\uc}{{\underline c}}

\newcommand{\uaa}{{\underline a}}
\newcommand{\ue}{{\underline e}}

\newcommand{\bb}{\begin{eqnarray}}
\newcommand{\be}{\end{eqnarray}}

\newtheorem{theorem}{Theorem}
\newtheorem{definition}{Definition}

\begin{document}

\title{Computational aspects and applications of a new
transform for solving the complex exponentials approximation
problem}
\author{Piero Barone
\thanks{ Istituto per le Applicazioni del Calcolo ''M. Picone'',
C.N.R., Viale del Policlinico 137, 00161 Rome, Italy; e-mail:
barone@iac.rm.cnr.it; fax: 39-6-4404306}}
\date{}
\maketitle

\section*{Abstract}

Many real life problems can be reduced to the solution of a complex
exponentials approximation problem which is usually ill posed.
Recently a new transform for solving this problem, formulated as a
specific moments problem in the plane, has been proposed in a
theoretical framework. In this work some computational issues are
addressed to make this new tool useful in practice. An algorithm is
developed and used to solve a Nuclear Magnetic Resonance
spectrometry problem, two time series interpolation and
extrapolation problems and a shape from moments problem.

{\it Key words:}
 complex moments problem; Pade'  approximants; logarithmic potentials;
 random determinants; pencils of matrices
{\it AMS classification:} 62M15, 30Exx

\pagestyle{myheadings}
\thispagestyle{plain}\markboth{P.Barone}{Complex exponentials
approximation}

\newpage

\section*{Introduction}

Many signal processing problems (see e.g. \cite{scharf}) can be
formulated as  a complex exponential interpolation problem (CEIP):
given the complex numbers $s_k,\; k=0,1,2,\dots 2p-1,$ to find
complex numbers $\{c_j,\xi_j\},j=1,\dots,p$ such that
\begin{eqnarray}s_k=\sum_{j=1}^{p}
c_j\xi_j^k,\;\;k=0,1,\dots,2p-1.\label{cei}\end{eqnarray} or,
equivalently \cite{hen2}, to find poles $\xi_j$ and corresponding
residues $r_j=c_j/\xi_j$ of the rational function $s(z)$ whose first
$2p$ Taylor coefficients at $z=0$ are $s_k,\; k=0,1,2,\dots 2p-1.$
The problem can be restated as a generalized eigenvalue problem as
follows. Let us consider Hankel matrices $U_0(\us)$ and $U_1(\us)$
given by
$$U_0(\us)=\left[\begin{array}{llll}
s_0 & s_{1} &\dots &s_{p-1} \\
s_{1} & s_{2} &\dots &s_{p} \\
. & . &\dots &. \\
s_{p-1} & s_{p} &\dots &s_{2p-2}
  \end{array}\right]$$
$$U_1(\us)=\left[\begin{array}{llll}
s_1 & s_{2} &\dots &s_{p} \\
s_{2} & s_{3} &\dots &s_{p+1} \\
. & . &\dots &. \\
s_{p} & s_{p+1} &\dots &s_{2p-1}
  \end{array}\right]$$
  where $\us=[s_0,\dots,s_{2p-1}].$
Because of (\ref{cei}), the following factorizations hold
$$U_0(\us)=VCV^T,\;\;
U_1(\us)=VCZV^T$$ where $V$ is the Vandermonde matrix based on
$(\xi_1,\dots,\xi_p)$,
$$C=diag\{c_1,\dots,c_p\} \mbox{ and }
Z=diag\{\xi_1,\dots,\xi_p\}.$$ Therefore $(\xi_j,\;j=1,\dots,p)$ are
the generalized eigenvalues of the pencil $ P=[U_1(\us),U_0(\us)]$
and $(c_j,\;j=1,\dots,p)$ are related to the generalized eigenvector
$\uu_j=V^{-T}\ue_j$ of $P$ by $c_j=\uu_j^T[s_0,\dots,s_{p-1}]^T$,
where $\ue_j$ is the $j-$th column of the identity matrix $I_{p}$ of
order $p$. A further equivalent formulation is based on  the complex
measure
$$S(z)=\sum_{j=1}^{p} c_j\delta(z-\xi_j),\;\;z, \xi_j\in
D,$$ where $D$ is a compact subset of $\C$ and $\delta$ is the Dirac
distribution. It turns out that for $k=0,1,2,\dots$
$$s_k=\int_Dz^kS(z)dz=\int\!\!\int_D(x+iy)^k S(x+iy)dx dy.$$
Therefore $s_k$ is the $k$-th  harmonic moment of the measure
$S$ and the complex exponential interpolation problem is equivalent
to a specific moment problem in the plane consisting in retrieving
the distribution $S$ from $s_k, \; k=0,1,2,\dots 2p-1$. Conditions
for existence and unicity of the solution are $det U_0(\us)\ne 0,
det U_1(\us) \ne 0$ (see e.g. \cite[Th.7.2c]{hen2}).

\noindent More realistically, by denoting in bold all random
quantities, let us consider the discrete stochastic process defined
by
\begin{eqnarray}{\bf a}_k=s_k+\mbox{\boldmath $\nu$}_k,\; k=0,1,2,\dots,n-1
\label{modale}\end{eqnarray} where $n\ge 2p$ and
$\mbox{\boldmath$\nu$}_k$ is a complex Gaussian zero-mean white
noise discrete process with known variance $\sigma^2$. We want
therefore to solve the complex exponential approximation problem
(CEAP) consisting in estimating $p$ and $\{c_j,\xi_j\},j=1,\dots,p,$
from a realization $a_k,\; k=0,1,2,\dots,n-1$ of ${\bf a}_k$. This
is equivalent, when $p$ is known, to solve a Pade' approximation
problem i.e. to compute the $[p,p-1]$ Pade' approximant of the
formal power series $f(z)=\sum_k a_kz^{-k}$, or to solving a
generalized eigenvalue problem for nonsquare pencils \cite{lem,but}
or a specific noisy moments problem in the plane. Even if $p$ were
known the problem would be quite difficult and usually ill-posed. A
wide literature exists on the subject. We can summarize some well
known facts as follows (see e.g. \cite{osb,gmv,donoho}). The problem
is optimally conditioned when $\xi_j$ are equispaced on the unit
circle. In this case in fact model (\ref{cei}) reduces to the
Fourier model which is an orthogonal one. Clusters of $\xi_j$ are
more difficult to estimate than well separated ones. Complex
exponentials with relatively small $|c_j|$ are more difficult to
estimate than those with relatively large weight.

Recently a new approach for solving the complex exponential
approximation problem in a stochastic framework was proposed
\cite{pb06}, which exploits the relation with generalized eigenvalue
problems and with moments problems outlined above but without
assuming to know $p$. It makes use of tools from the theory of
logarithmic potential with external fields \cite{saff} and the
theory of random polynomials \cite{randpol,ham} and provides an
estimate of $p$ and point and interval estimates of $\{c_j,\xi_j\}.$

In this work some computational and numerical issues are addressed
to make this new tool useful in practice. An algorithm is developed
and tested on well known difficult problems.

 \noindent The paper is organized as follows.
In Section 1 the method introduced in \cite{pb06} is shortly
summarized. In Section 2 the proposed algorithm  is discussed.  In
Section 3 the algorithm is used to solve a Nuclear Magnetic
Resonance spectrometry problem, a time series interpolation and
extrapolation problem and a shape from moments problem providing
some comparisons with existing methods.

\section{The new transform}

Starting from  ${\bf a}_k,\; k=0,1,2,\dots,n-1$, assuming $n$ even,
let us consider the stochastic CEIP (i.e. a CEIP for each
realization of $\{{\bf a}_k\}$)
\begin{eqnarray}{\bf a}_k=\sum_{j=1}^{n/2}
{\bf c}_j\mbox{\boldmath
$\xi$}_j^k,\;\;k=0,1,\dots,n-1.\label{ceip}.\end{eqnarray} and the
associated random measure
\begin{eqnarray}{\bf S}_n(z,\sigma)=\sum_{j=1}^{n/2}{\bf c}_j\delta(z-\mbox{\boldmath $\xi$}_j)
\label{ranmes}.\end{eqnarray} Let us also define  the random Hankel
$\frac{n}{2}\times \frac{n}{2}$ matrices $U_0({\bf
\uaa}),\;\;U_1({\bf \uaa})$, where ${\bf \uaa}=[{\bf a}_0,\dots,{\bf
a}_{n-1}]$. The generalized eigenvalues $\mbox{\boldmath
$\xi$}_j,\;j=1,\dots,n/2$ of the random pencil ${\bf P}=[U_1({\bf
\uaa}),U_0({\bf \uaa})]$   satisfy the equation
$${\bf p}_{n/2}(z)=det[U_1({\bf
\uaa})-zU_0({\bf \uaa})]=0$$ where ${\bf p}_{n/2}(z)$ is a random
polynomial. We can then consider the expected value of the (random)
normalized counting measure on the zeros $\mbox{\boldmath
$\xi$}_j,\;j=1,\dots,n/2$ of this polynomial (condensed density,
\cite{ham,randpol}):
$$h_n(z)= \frac{2}{n}E\left[\sum_{j=1}^{n/2}\delta(z-\mbox{\boldmath
$\xi$}_j)\right].
$$
In  \cite{barja} it was proved that, when $\us=\underline{0}$, in
the limit for $n\rightarrow \infty$ the condensed density is a
distribution supported on the unit circle and it can be proved
(\cite{pb06}) that in the limit for $\sigma\rightarrow \infty$ the
generalized eigenvalues $\mbox{\boldmath $\xi$}_j$ tend to
concentrate on the unit circle and, in the limit for
$\sigma\rightarrow 0$, they concentrate around  the true
$\xi_j,j=1,\dots,p$. It is therefore evident that in order to solve
CEAP, the first issue to address is the identifiability one. If the
Signal-to Noise  ratio (SNR) is not large enough with respect to the
signal structure as discussed in the introduction, there is no hope
to solve CEAP. The first step of the method introduced in
\cite{pb06}  provides a tool for assessing if CEAP is solvable based
on the properties of the condensed density of the generalized
eigenvalues $h_n(z)$. More precisely we give the following:
\begin{definition}
The measure $S(z)$  is identifiable from ${\bf a}_k,k=0,\dots,n-1$
if $\exists \;\;r_k>0,k=1,\dots,p$ such that
\begin{itemize}
\item
$h_n(z) \mbox{ is unimodal in }  N_k=\{z\; \| \;\;|z-\xi_k|\le
r_k\}$
\item
$ \bigcap_{k=1}^{p}N_k=\emptyset$
\end{itemize}
\end{definition}
\noindent The following result, proved in \cite{pb06}, gives the
relation between ${\bf S}_n(z,\sigma)$, and the unknown measure
$S(z)$
\begin{theorem}
If $S(z)$ is identifiable from $\uaa$ then
$$\int_{N_h(r_h)}E[{\bf S}_n(z,\sigma)]dz= c_h+o(\sigma),\;\;h=1,\dots,p$$
and $$\int_{A}E[{\bf S}_n(z,\sigma)]dz= o(\sigma),\;\;\forall
A\subset D-\bigcup_jN_j(r_j).$$\label{teo33}
\end{theorem}
As in the limit for $\sigma\rightarrow 0$, the condensed density
tends continuously to a distribution supported on the true
$\xi_j,j=1,\dots,p$, it does exist $\sigma$ small enough to make
$S(z)$ identifiable from ${\bf \uaa}$ and in this case we can use
the random measure ${\bf S}_n(z,\sigma)$ to estimate $S(z)$ by using
Theorem \ref{teo33}. To perform this program we need two steps. The
first one consists in either to check the identifiability of the
measure $S(z)$ from ${\bf \uaa}$  or to properly design the
experiment (i.e. to choose $n$ and $\sigma$) in order to get
identifiability. The second step consists in building an estimator
of ${\bf S}_n(z,\sigma)$.

About the first step we notice that of course the function $h_n(z)$
cannot be computed because we do not know  $\us$ i.e. the mean of
${\bf \uaa}$.  However, assuming to know $\us$, we can use $h_n(z)$
to state whether $S(z)$ is identifiable from the data. Unfortunately
even in the Gaussian assumption the analytic computation of $h_n(z)$
is hard. However it can be approximated (\cite{pb06}) by
$$\tilde{h}_n(z)=\frac{1}{2\pi n}\Delta\sum_{\mu_j(z)> 0
}\log(\mu_j(z))$$  where $\Delta$ is the Laplacian operator acting
on $z$ and $\mu_j(z)$ are the eigenvalues of
\begin{eqnarray}(U_1(\us)-zU_0(\us))\overline{(U_1(\us)-zU_0(\us))}+
\frac{n\sigma^2}{2}A(z,\overline{z})\label{appr}\end{eqnarray} where
$\us=[s_0,\dots,s_{n-1}]$,
$$A(z,\overline{z})=\left[
\begin{array}{lllll} 1+|z|^2 & \;\;-z &\; 0&\dots&0
\\ -\overline{z}\;\; & 1+|z|^2& \;\;-z &\;\;0&\dots\\ . & . & . & .&.\\
  0 &\dots&0& -\overline{z}\;\; & 1+|z|^2 \end{array}\right]\in
  \C^{\frac{n}{2}\times\frac{n}{2}}$$ and overline denotes
  conjugation.

\noindent\underline{Remark.} From  equation (\ref{appr}) it follows
that $n$ should not be as large as possible to get the best
estimates of $S(z)$. In fact too many data will convey too much
noise which could mask the signal $s_k$.

\noindent We
have therefore a tool either to check identifiability or to design
properly the experiment. In most real problems we have some prior
information about the unknown measure $S(z)$. We can then compute
$\tilde{h}_n(z)$ for several candidate measures  compatible with our
prior information and choose values $n$ and $\sigma$ that make the
candidate measures identifiable.

We now move to the second step of the procedure consisting in
estimating the random measure ${\bf S}_n(z,\sigma)$ and extracting
from it the required information. If we have $R$ samples  from the
data discrete stochastic process ${\bf \uaa}$ we can estimate
$E[{\bf S}_n(z,\sigma)]$ by solving CEIP for each sample
$\uaa^{(r)},\;r=1,\dots,R,$ i.e. finding
$(c_j^{(r)},\xi_j^{(r)}),j=1,\dots,n/2,$ such that
$a^{(r)}_k=\sum_{j=1}^{n/2}
c_j^{(r)}(\xi_j^{(r)})^k,\;\;k=0,1,\dots,n-1$ and then taking the
sample mean
$$\frac{1}{R}\sum_{r=1}^R\sum_{j=1}^{n/2}c_j^{(r)}\delta(z-\xi_j^{(r)}).
$$  If only one sample  is available we can use the following method proposed in
\cite{pb06}.
 We notice first that in order to cope with the Dirac distribution
  appearing in the definition of ${\bf S}_n(z,\sigma)$, it is convenient to use
  an alternative expression given by (see \cite{barja})
$${\bf S}_n(z,\sigma)=\frac{1}{4\pi}
\Delta\sum_{j=1}^{n/2}{\bf c}_j\log(|z-\mbox{\boldmath $\xi$}_j|^2)
.$$  Then we build
 independent replications of the data process
(pseudosamples) by defining
$${\bf a}_k^{(r)}={\bf a}_k+\mbox{\boldmath
$\nu$}_k^{(r)},\;\;k=0,\dots,n-1;\;\;\;r=1,\dots,R$$ where
$\{\mbox{\boldmath $\nu$}_k^{(r)}\}$ are i.i.d. zero mean complex
Gaussian variables with variance $\sigma'^2$ and therefore ${\bf
a}_k^{(r)}$ have variance $\tilde{\sigma}^2=\sigma^2+\sigma'^2$. We
then define the estimator, conditioned to ${\bf \uaa}$
\begin{eqnarray}
{\bf S}_{n,R}^c(z,\tilde{\sigma})=\frac{1}{2\pi
R}\Delta\left(\sum_{r=1}^R\sum_{j=1}^{n/2} {\bf
c}_j^{(r)}\log(|z-\mbox{\boldmath $\xi$}_j^{(r)}|)
\right)\label{ptransf}
\end{eqnarray} where
$({\bf c}_j^{(r)},\mbox{\boldmath $\xi$}_j^{(r)}),j=1,\dots,n/2$ are
the solution of CEIP for the pseudodata ${\bf
a}_k^{(r)},\;\;k=0,\dots,n-1$ which are computable by a MonteCarlo
procedure given $\uaa$. In \cite{pb06} the following theorem is
proved
\begin{theorem}
Let $M(z)$ and $M_c(z)$ be the mean squared error of ${\bf
S}_{n}(z,\sigma)$ and ${\bf S}_{n,R}^c(z,\tilde{\sigma})$
respectively. In the limit for $\sigma\rightarrow 0$, it exists
$\sigma'$ and $R(\sigma')$ such that $\forall R\ge R(\sigma')$,
$M_c(z)<M(z)\;\; \forall z$.
\end{theorem}

\noindent In order to estimate  $(c_j,\xi_j),\;j=1,\dots,p,$ we make
use of Theorem \ref{teo33}. In fact, if $S(z)$ is identifiable,
there exist disjoint sets $N_k,k=1,\dots,p$ such that $\left|
\int_{N_k}E[{\bf S}_n(z,\sigma)]dz\right|\gg\sigma$, and each of
them should include one and only one $\xi_k$. Therefore looking at
the sets $A$ such that $\left| \int_{A}{\bf
S}_{n,R}(z,\tilde{\sigma}))dz\right|\gg\sigma$ it is possible to
identify $\hat{p}$ disjoint sets $\hat{N}_k$ which possibly include
the true $\xi_k$. This can be done by computing  a discrete
transform  by
 evaluating ${\bf S}_{n,R}^c(z,\tilde{\sigma})$ on a suitable
lattice $L=\{(x_i,y_i),i=1,\dots,N\}$ by taking a discretization of
the Laplacian operator, giving rise to a matrix
$\mathsf{P}(\tilde{\sigma})\in\Re_+^{(N\times N)}$ - the
$\mathsf{P}$-transform of the vector $\uaa$ -  such that
${\mathsf{P}(i,j,\tilde{\sigma})}={\bf S}_{n,R}^c(x_i+iy_j)$. The
$\hat{p}$ relative maxima of the absolute value of the
$\mathsf{P}$-transform are then computed as well as disjoint
neighbors $\hat{N}_k$ centered on them. Estimates
$(\hat{c}_k,\hat{\xi}_k)$ of $(c_k,\xi_k)$ are obtained by averaging
the $(c^{(r)}_j,\xi^{(r)}_j)$ which belong to each $\hat{N}_k$. The
name "transform" is justified by observing that to the vector $\uaa$
we associate the matrix $\mathsf{P}$ (direct transform), and to the
matrix $\mathsf{P}$ we associate the vector whose components are
$$\hat{a}_k=\sum_{j=1}^{\hat{p}}\hat{c}_j\hat{\xi}_j^k\rightarrow a_k,\mbox{ when }\sigma\rightarrow
0 $$ (inverse transform).

\section{The algorithm}

 The method for estimating  the unknown parameters $p$, $\{(c_j,\xi_j),
,j=1,\dots,p\}$ outlined in the previous section is quite expensive
and delicate from the numerical point of view. In this section we
discuss the main issues to be addressed to implement the basic
method and  suggest a new approach which mimics the basic one giving
rise to a fast and reliable algorithm.

 The
computation of the $\mathsf{P}$-transform is the most critical part
of the whole procedure. There are many algorithms to compute
$(\hat{c}_j,\hat{\xi}_j)$ based on different approaches (see e.g.
\cite{elad,bama2} for  short reviews) which are useful in different
applied contexts. If  computational burden is the principal issue
and the geometric structure of the unknown measure $S(z)$ is simple,
extremely fast algorithms based on the generalized orthogonality of
Pade' polynomials can be used to compute $(\hat{c}_j,\hat{\xi}_j)$
(\cite[pg.631-632]{hen2},\cite{brez}). If clusters of poles can be
expected it is better to solve the generalized eigenvalue problem
e.g. as discussed in \cite{osb} and \cite{elad} where several
advanced methods are presented or \cite{gmv}, where the Hankel
structure of the pencil $P=[U_1(\uaa),U_0(\uaa)]$ is taken into
account to speed up the computation and QR factorization and QZ
iteration are used as well as a suitable diagonal scaling of the
pencil $P$, for achieving numerical stability. An even more
expensive method is described in \cite{schu} where a total least
squares approach is used taking into account the Hankel structure
and the noise affecting the elements of $P$. A classical approach is
given by Prony's method \cite{pr} which splits the problem in three
parts by solving two linear least squares problems with Toeplitz and
Vandermonde structure respectively and a polynomial rooting problem.
Fast codes for all these sub-steps do exist \cite[sect.4.6,4.7]{gl}
as well as total least squares \cite{huf} and structured total least
squares algorithms \cite{lem}.

\noindent A further complication is due to the fact that for
computing the $\mathsf{P}$-transform  $R$ generalized eigenvalue
problems have to be solved. An effective compromise between accuracy
and speed of computation is given by the following procedure:
\begin{itemize}
\item
 compute $(c_j^{(0)},\xi_j^{(0)}),\;j=1,\dots,n/2,$ by solving the generalized
eigenvalue problem for the  pencil $P$ by one of the accurate
methods quoted above. If the method described in \cite{gmv} is used
the computational cost of this step is $O((\frac{n}{2})^3)$
\item select the generalized eigenvalues $\xi_j^{(0)},\;j=1,\dots,\tilde{p}$
corresponding to the
$\tilde{p}$ largest values  $|c_j^{(0)}|$ where  $\tilde{p}$ is an
upper bound of  $p,\;\;\tilde{p}\le\frac{n}{2}$
\item for each pseudosample $\uaa^{(r)}$  compute
the coefficients of the polynomial
$$p(z)=\det[U_1(\uaa^{(r)})-zU_0(\uaa^{(r)})]$$ by the first step of Prony's
method. This requires $O((\frac{n}{2})^2)$ flops because of the
Hankel structure of $U_0(\uaa^{(r)})$
\item
to compute $\xi_j^{(r)},\;j=1,\dots\tilde{p}$, apply a fast
iteration such as e.g. Laguerre method \cite{wil} to the polynomial
$p(z)$, taking as initial values
$\xi_j^{(0)},\;j=1,\dots,\tilde{p}$. Usually it converges in few
iterations, therefore it costs $O(\frac{n}{2}\tilde{p})$ flops or
less if the Horner scheme to compute the polynomial derivatives is
implemented through the fast Fourier transform \cite{fact}.
\item to compute $c_j^{(r)},\;j=1,\dots\tilde{p}$,
apply the third step of Prony's method, forming the  Vandermonde
matrix of $\xi_j^{(r)},\;j=1,\dots\tilde{p}$ and solving a least
squares problem e.g. by  LSQR method \cite{ps}, which is a good
compromise between accuracy and computational speed. Usually a few
iterations are sufficient, therefore it costs
$O(\frac{n}{2}\tilde{p})$ flops.

\item The last step for computing the $\mathsf{P}$-transform
 consists in evaluating the summation in (\ref{ptransf}) and then
 computing a discrete Laplacian. This can be the most
expensive part of the procedure because the summation must be
computed for each $z_i=(x_i,y_i),\;i=1,\dots,N$ of the lattice $L$.
Therefore we need $O(R\tilde{p}N^2)$ flops. However we notice that
we only need an estimate of the local maxima of the absolute value
of the $\mathsf{P}$-transform. These are likely to be close to the
centroids of poles clusters and their value is a monotonic
increasing function of the corresponding $|c_j^{(r)}|$. A fast way
to estimate them consists then in applying a clustering method, such
e.g. k-means, to the $R\tilde{p}$ vectors of $\R^3$
$$[\Re{\xi_j^{(r)}},\Im{\xi_j^{(r)}},
|c_j^{(r)}|],r=1,\dots,R,\;j=1,\dots,\tilde{p}$$ looking for
$\tilde{p}$ clusters. The clustering algorithm can be initialized by
$$[\Re{\xi_j^{(0)}},\Im{\xi_j^{(0)}},
|c_j^{(0)}|],\;j=1,\dots,\tilde{p}$$ computed in the first two
steps. We then compute
\begin{eqnarray}
S_{n,R}^c(z,\tilde{\sigma})=\frac{1}{2\pi
R}\Delta\left(\sum_{r=1}^R\sum_{j=1}^{\tilde{p}}
c_j^{(r)}\log(|z-\xi_j^{(r)}|)   \right)
\end{eqnarray}
for $z\in N_k,\;\;k=1,\dots,\tilde{p}$ where $N_k$ is a small
regular mesh of points with size $\delta$, centered on the centroid
of the $k-$th cluster. Finally, after theorem \ref{teo33}, we select
the $\hat{p}\le\tilde{p}$ clusters such that
$$\left|\sum_{z_h\in
N_k}S_{n,R}^c(z_h,\tilde{\sigma})\delta^2\right|>K\sigma$$ where
$K>1$. Estimates $(\hat{c}_j,\hat{\xi}_j),\;j=1,\dots,\hat{p}$ of
$(c_j,\xi_j),\;j=1,\dots,p$ are then obtained by averaging the
$(c_j^{(r)},\xi_j^{(r)})$ which belong to the selected clusters. The
computational cost of the clustering algorithm and the computation
of $S_{n,R}^c(z,\tilde{\sigma})$ is $O(R(\tilde{p})^2)$ flops.

\end{itemize}

\noindent Summing up we can solve the CEAP problem in
$O((\frac{n}{2})^3)+O(R(\frac{n}{2})^2)$ flops. In most applications
$R<n$ is enough to get good results, therefore $O(n^3)$ flops is a
reasonable upper bound for solving the problem in most cases (fast
method). In a few particularly difficult problems the computation of
$\xi_j^{(r)},\;j=1,\dots\tilde{p}, r=1,\dots,R$ is better performed
by the same accurate methods used for $r=0$. In these cases the
computational burden becomes $O(n^4)$ (slow method).

\section{Numerical experiments}

In order to appreciate the behavior of the proposed algorithm in
practice, four examples on real and synthetic data are presented.
 The first one copes with the classic problem of
quantification of Nuclear Magnetic Resonance spectra (see e.g.
\cite{vin1,vin2}) which is usually solved by ad hoc methods
requiring visual inspection by the operator.
 The second example is an interpolation-extrapolation
problem on a synthetic time series used in the 2004 Competition on
Artificial Time Series, organized in the framework of the European
Neural Network Society \cite{cats}. Comparisons with the results
obtained by  participants are provided.
 The third example is an interpolation problem of a real
acoustic signal with a missing fragment. The aim here is to
reconstruct the missing part in order to make the reconstructed
signal to sound as the original. The last example is a shape from
moments reconstruction problem. It turns out that the identification
of a polygonal region in the plane from its complex moments can be
formulated as a specific CEIP \cite{davis,gmv}. Synthetic data sets
are generated and the results are compared with those obtained in
\cite{elad} when the number of the polygon vertices is known.
Moreover the case when the number of vertices is unknown is also
addressed.

We notice that several hyperparameters have  to be chosen e.g. the
upper bound $\tilde{p}$ of $p$,  the number $R$ of pseudosamples,
the variance $\sigma'^2$ of $\{\mbox{\boldmath $\nu$}_k^{(r)}\}$ and
the constant $K$. Moreover one of the most critical hyperparameter
is the number $n$ of data points, as noted in the Remark in section
2. Usually we can only cut some data in order to reduce the noise.
In order to select good hyperparameters a performance criterion is
chosen and the method is applied for several values of the
hyperparameters in suitable intervals. Then those that give the best
value of the performance criterion are used to compute the final
results. The performance criterion is problem dependent. However a
standard residual analysis provides usually a good basis to build up
a good criterion. In the following the number of residuals whose
absolute value is larger than $\sigma$ is used as performance
criterion.

\subsection{NMR spectroscopy}

\begin{figure}
\begin{center}
\hspace{.2cm}\epsfig{file=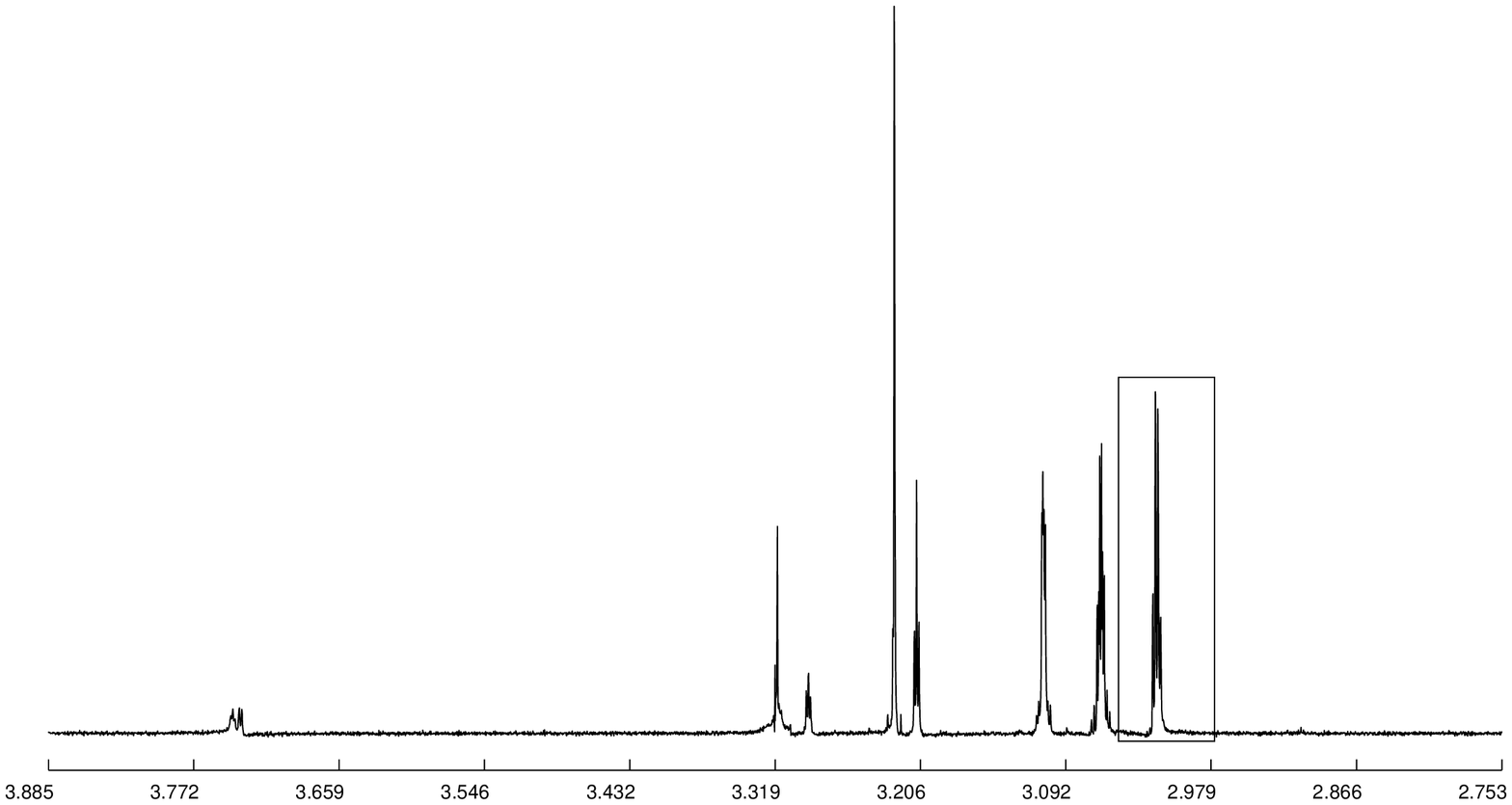,width=5cm,height=14cm,angle=270}
\hspace{.2cm}\epsfig{file=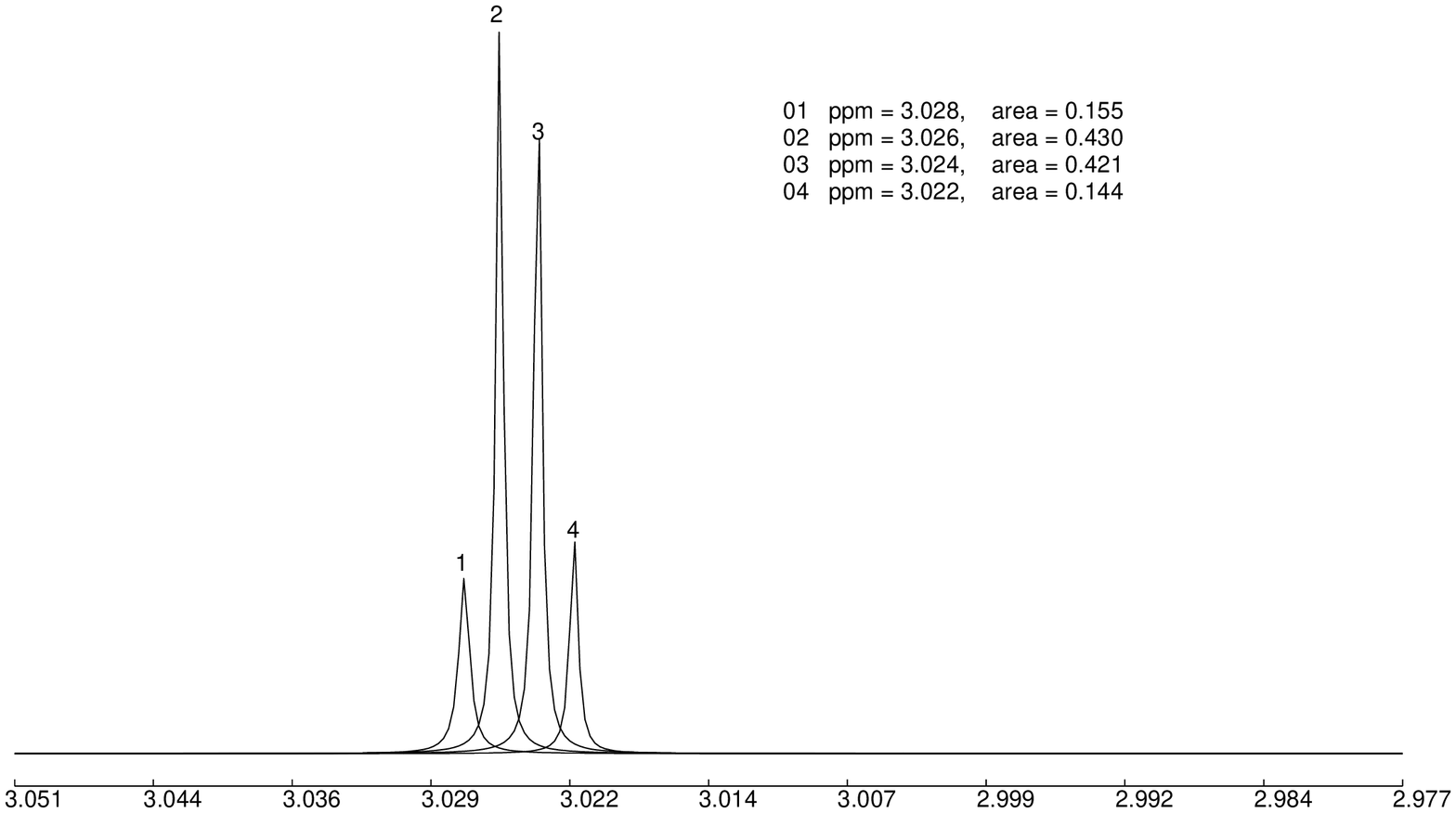,width=5cm,height=14cm,angle=270}
\vspace{.2cm}
\hspace{.2cm}\epsfig{file=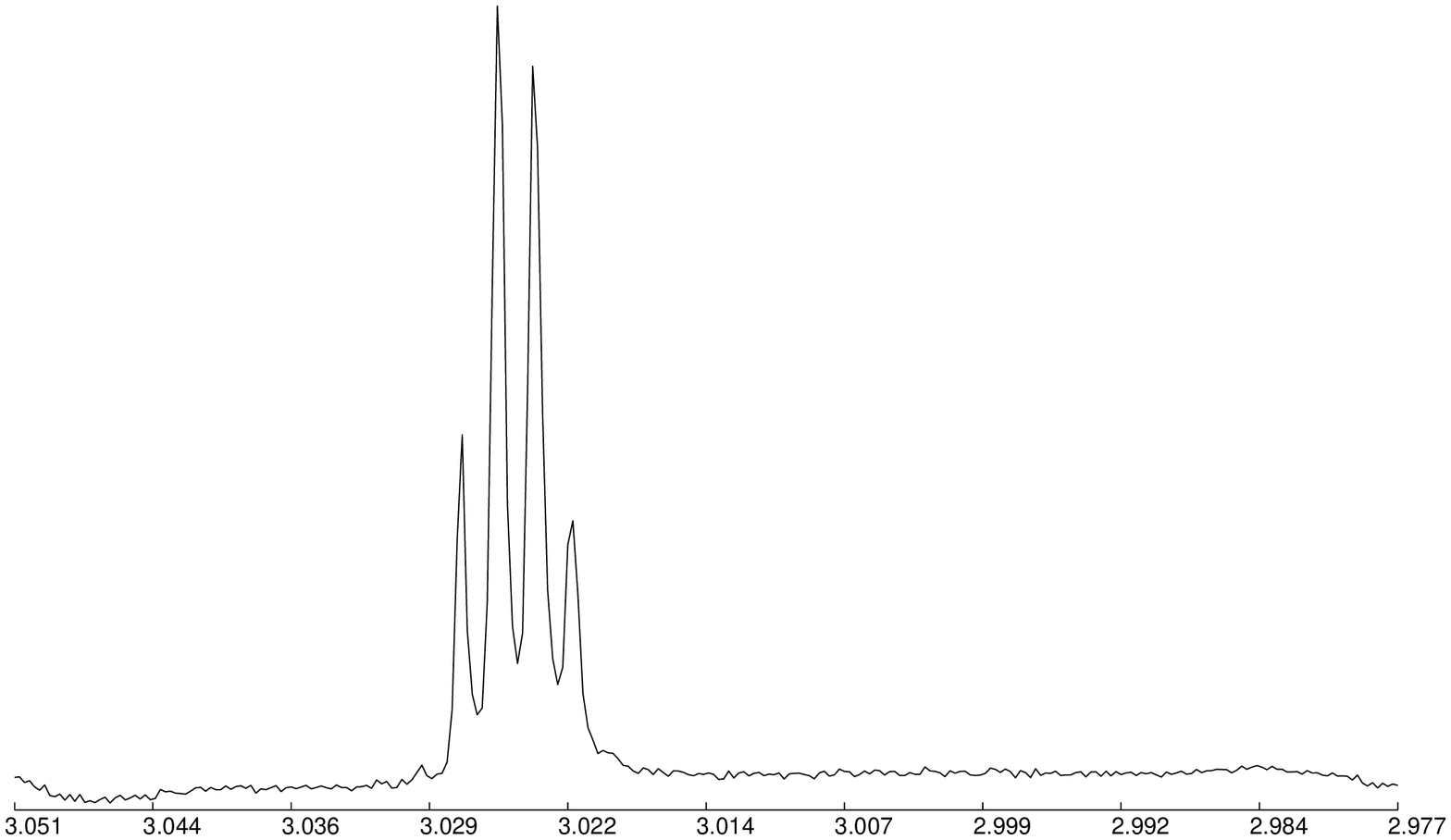,width=5cm,height=14cm,angle=270}
\hspace{.2cm}\epsfig{file=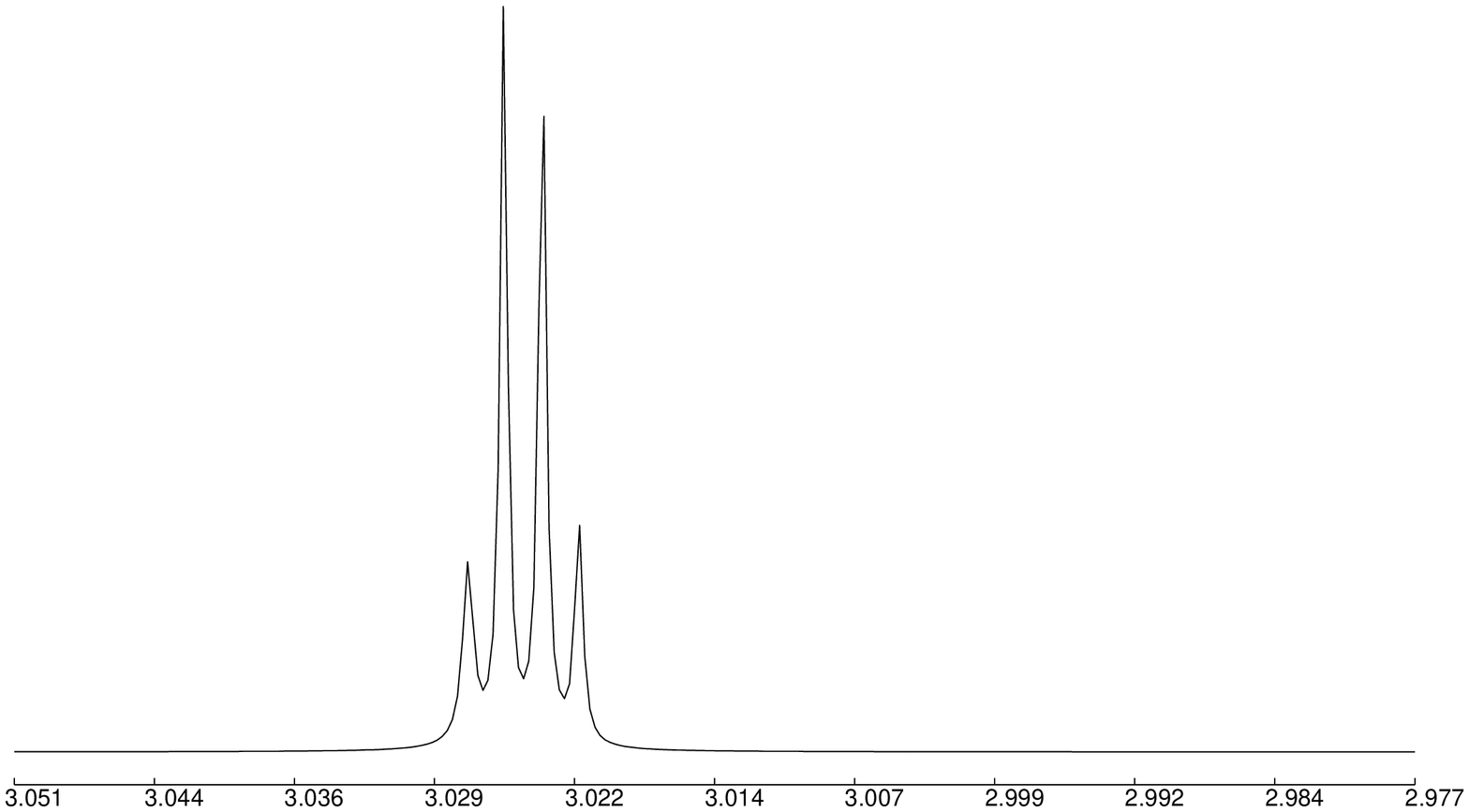,width=5cm,height=14cm,angle=270}
\caption{Quantification of NMR spectra. Top: the NMR Fourier
spectrum and the ROI. Middle-top: the estimated Lorentian lines in
the ROI and their areas. Middle-bottom: the Fourier spectrum in the
ROI. Bottom: the sum of the estimated Lorentian lines.} \label{fig1}
\end{center}
\end{figure}

In the top part of fig.\ref{fig1}  a $^1H$  Magnetic Resonance
absorption spectrum from a diluted aqueous solution of the
tripeptide Glutathione in its reduced form (GSH) is shown. It is
computed by taking the discrete Fourier transform (DFT) of a Free
Induction Decay (FID) signal of $4096$ data points. In ideal
conditions the physical model for the FID is a linear combination
with  positive weights of complex exponentials. The absorption
spectrum is the real part of the Fourier transform of the FID. It
turns out that it is given by a linear combination of Lorentian
functions. The spectroscopist is interested in estimating the
parameters characterizing these Lorentian lines, namely their modes,
widths and relative areas which are simply related respectively to
the argument of the complex exponentials modeling the FID, to their
absolute value and to weights associated to them. In a real
experimental setup the ideal conditions are no longer true. Standard
methods, implemented on most spectrometers,   fit each peak of the
absorption spectrum with a Lorentian function. If the peaks are
close,  a very ill conditioned non linear problem has to be solved
which can heavily depend on the interactive choices of the
spectroscopist to initialize the procedure. A better alternative is
provided by time-domain methods (see e.g.\cite{vin1,vin2}) which
exploit the fact that the FID can be modeled by complex
exponentials. The problem can still be very ill conditioned.
However, if the SNR is large enough, reasonable estimates of the
quantities of interest can be obtained by solving the CEAP problem
for the FID by the proposed method, which provides a global stable
solution and no longer requires critic interaction with the
spectroscopist.

The analysis is performed in the interval of the spectrum marked by
the rectangle in the top part of fig.\ref{fig1}. A quadruplet whose
areas are in the ratios $1:3:3:1$ is the theoretical reference. The
frequencies are measured in  parts per million (ppm). The standard
interactive procedure provides an estimate of the areas of the four
peaks such that their ratios are
 $1.02:3.21:3.15:1$. In order to apply the
proposed method the FID is first filtered by a pass-band Fourier
filter \cite{belk,bama2}. In the middle-bottom part of the same
figure, the absorption spectrum of the filtered FID is shown. When
the main peaks of the spectrum are clustered and the clusters are
well separated, it is in fact possible to split the analysis by
filtering out from the FID all the frequencies but those belonging
to a given interval \cite{neu}. The filtered FID  is given by only
$300$ data points and the proposed method was applied to solve the
CEAP for it. The results are shown on the middle-top part of
fig.\ref{fig1}. Four estimated Lorentian lines marked $1-4$ are
plotted and their areas are reported in the legend as well their
modes in ppm. The ratios of the areas are $1.07:2.98:2.92:1$ which
compare favorably with those estimated interactively. On the bottom
part of the figure the weighted sum of the four Lorentian lines is
plotted. The agreement with the zoomed absorption spectrum on the
middle-bottom  part is quite good.

\subsection{Time series interpolation and extrapolation}

In order to apply the proposed method to solve extrapolation
problems it is enough to solve a CEAP  for the measured data and
then evaluate the model on the extrapolation abscissas. To solve an
interpolation problem we notice that, in the noiseless case, we can
consider the segments of data before and after the missed segment as
produced by the same model (\ref{cei}) for a set of indices $k$
displaced by a fixed quantity $q$. It is easy to show that the
generalized eigenvalues and eigenvectors are invariant for such a
displacement. Therefore we can solve two separate CEAPs  for the
observed segments, and apply the proposed method to the pooled
generalized eigenvalues and eigenvectors. We need only to modify the
Vandermonde matrix for computing $c_j^{(r)}$ in the last step to
take into account the gap in the observations. Assuming that each
segment has $n$ observations we have $\uc^{(r)}=V^\dagger \uaa,$
where
$$V=\left[\begin{array}{llll}
 1& 1 &\dots &1 \\
\xi_1^{(r)} & \xi_2^{(r)} &\dots &\xi_{\tilde{p}}^{(r)} \\
. & . &\dots &. \\
. & . &\dots &. \\
(\xi_1^{(r)})^{n-1} & (\xi_2^{(r)})^{n-1} &\dots &(\xi_{\tilde{p}}^{(r)})^{n-1} \\

(\xi_1^{(r)})^{n+q-1} & (\xi_2^{(r)})^{n+q-1} &\dots &(\xi_{\tilde{p}}^{(r)})^{n+q-1} \\
. & . &\dots &. \\. & . &\dots &. \\
(\xi_1^{(r)})^{2n+q-1} & (\xi_2^{(r)})^{2n+q-1} &\dots
&(\xi_{\tilde{p}}^{(r)})^{2n+q-1}
  \end{array}\right]$$
 The interpolated values are then obtained by
 $$\uaa_{int}=V\hat{\uc},\;\;\;V=\left[\begin{array}{llll}
\hat{\xi}_1^{n} & \hat{\xi}_2^{n} &\dots &\hat{\xi}_{\hat{p}}^{n} \\
. & . &\dots &. \\
. & . &\dots &. \\
\hat{\xi}_1^{n+q-1} & \hat{\xi}_2^{n+q-1} &\dots
&\hat{\xi}_{\hat{p}}^{n+q-1}
  \end{array}\right].$$

\begin{figure}
\center{\fbox{\epsfig{file=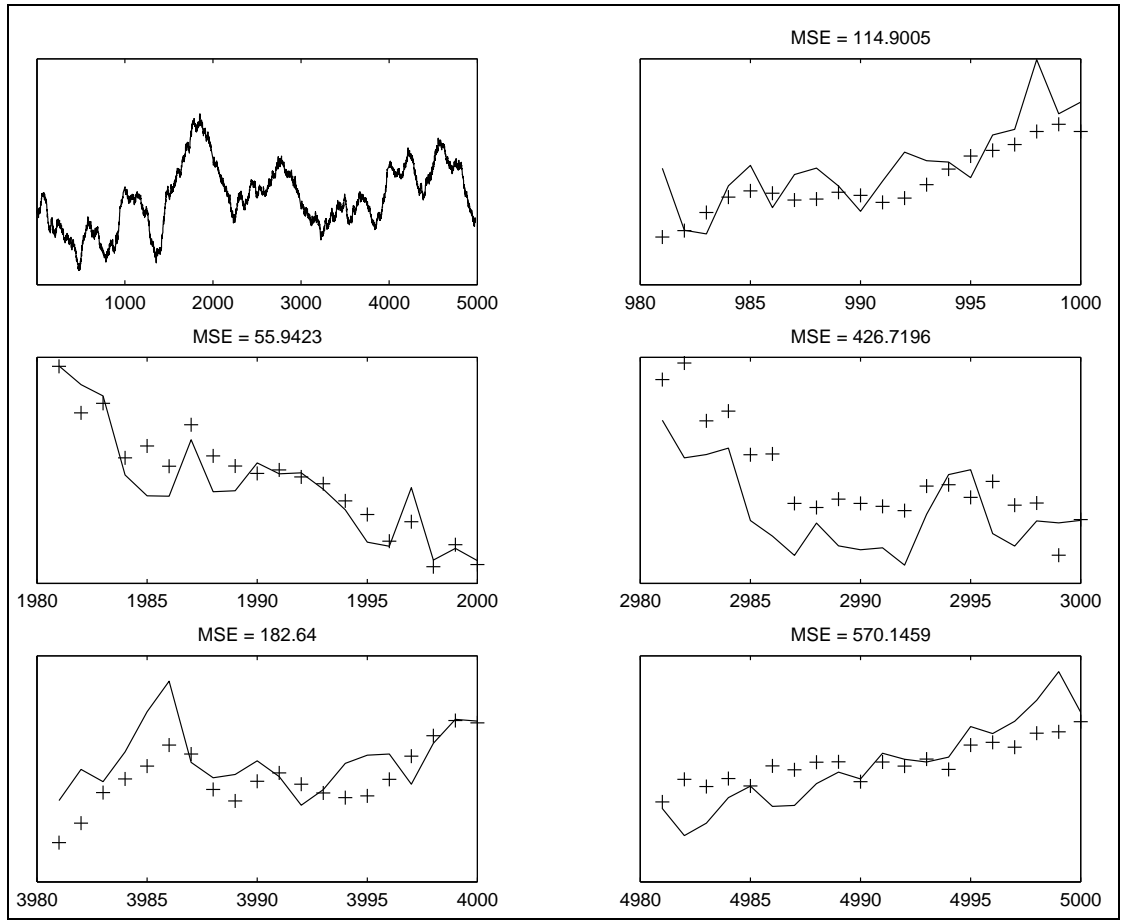,height=4.7in}}} \caption{Top
left: time series with five missing intervals. True values on each
interval (-); interpolated values (+). Total MSE on the first four
intervals = 195. Total MSE on the five intervals = 270. }
\label{fig2}
\end{figure}

\noindent The first example in this class of problems copes with a
time series of $5000$ samples with $100$ missing values at times
$981-1000, 1981-2000,2981-3000,3981-4000,4981-5000$. Therefore we
want to solve four interpolation and one extrapolation problems. As
the data are synthetic the truth is known and the results obtained
by $17$ methods are reported in \cite{cats} where the mean squared
error (MSE) for the interpolation problems and the interpolation +
extrapolation problems are reported. It can be argued that the MSE
is not the best discrepancy measure for this data set because a fit
with a smoothing cubic spline gives results better than all of the
$17$ quoted methods for the interpolation + extrapolation problems
and better than $15$ of them for the extrapolation problem.
Therefore we want to see how much the proposed method is able to
improve on the solution provided by the cubic spline. We then apply
the method to the residual obtained by subtracting the smoothing
spline from the data. In fig.\ref{fig2} top left the full time
series with missing data is plotted. The other plots show the true
values and the reconstructed ones on each missed data interval. The
$MSE_{100}=270$ and $MSE_{80}=195$ have to be compared with
$MSE_{100}=408$ and $MSE_{80}=222$ which are the best results
obtained in \cite{cats} by two different methods among the $17$
considered.

\begin{figure}
\center{\fbox{\epsfig{file=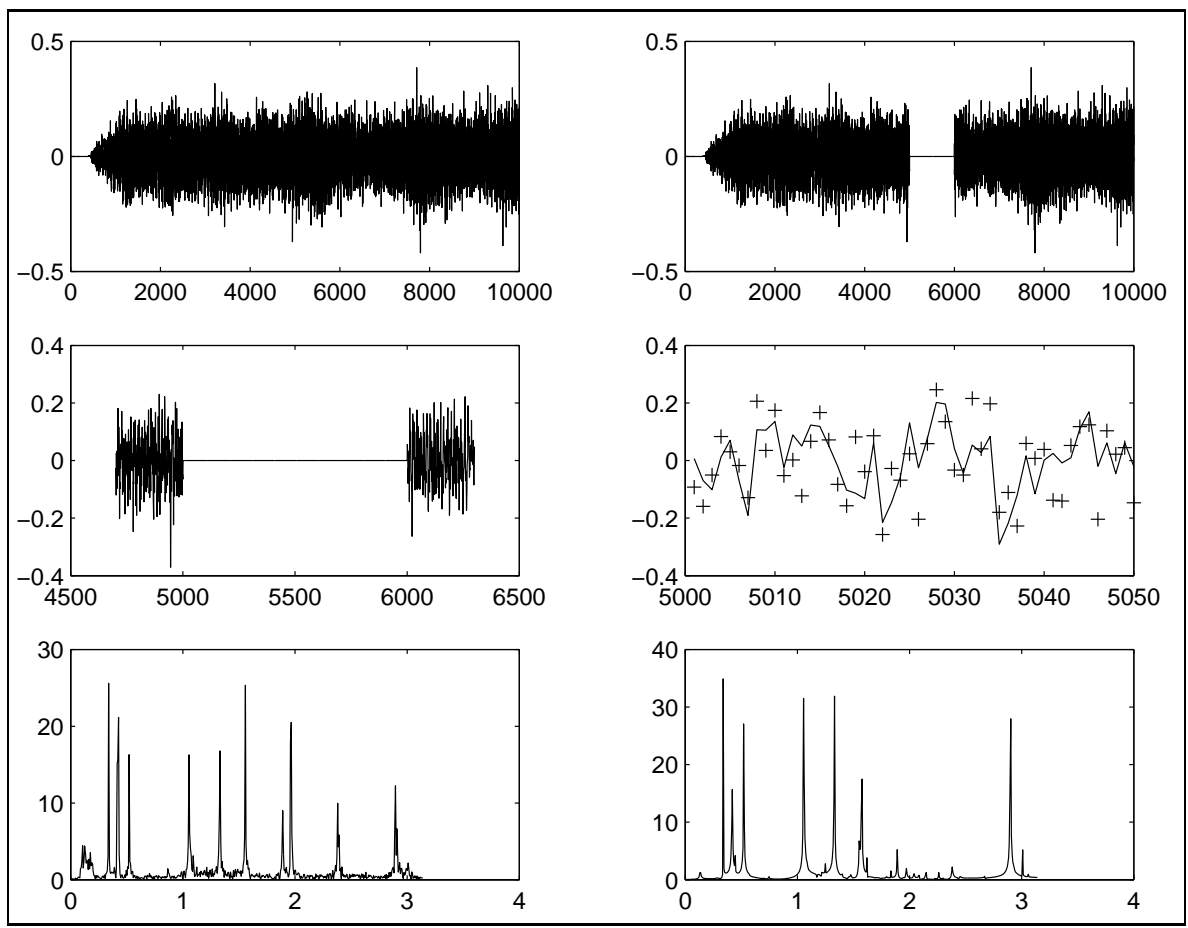,height=4.7in}}} \caption{ Top
left: a segment of an audio signal; top right: the part missed is
shown; middle left: the data used to interpolate; middle right: a
fraction of the interpolated data (+) are superimposed to the
unknown ones (-); bottom left: Fourier spectrum of the missed part;
bottom right: Fourier spectrum of the reconstructed data. }
\label{fig3}
\end{figure}

\noindent The second example is illustrated in fig. \ref{fig3}. An
audio signal, corresponding to a ringing bell, made up of $50000$
samples at $11025$ Hz is considered. The first $10000$ samples are
plotted in the top left part of the figure. A fragment of $1000$
samples are put to zero as shown in the top right part of the
figure. The method is applied to interpolate the missing fragment.
Two data sets made up of $300$ samples each before and after the
missing data are considered as shown in the middle left part of the
figure. The results are shown in the middle right part of the figure
where $50$ missed data are plotted superimposed to the interpolated
values. Even if the fit is not impressive most of the  main spectral
characteristics of the signal are well reproduced as shown in the
bottom part of the figure where the Fourier spectrum of the original
complete signal is plotted on the left, and the Fourier spectrum of
the complete signal with the missing fragment replaced by the
interpolated values is shown on the right. The sound produced by the
reconstructed signal is almost undistinguishable from the original
one.

\begin{table}[tbh]
\begin{scriptsize}
\begin{center}
\begin{tabular}{|c||c|c|c|c|c|}
\hline\hline
&$\sigma$&$RMSE \cite{elad}$&$RMSE$\\
\hline\hline
Star shape&1e-3&5.74e-2&3.68e-2\\
\hline
&1e-4&1.74e-2&1.02e-2\\
\hline
&1e-5&1.71e-3&1.05e-3\\
\hline \hline
C shape&1e-3&4.46e-3&4.30e-3\\
\hline
&1e-4&4.51e-4&4.27e-4\\
\hline
&1e-5&4.59e-5&4.28e-5\\
\hline \hline
\end{tabular}
\end{center}
\caption{For the star shaped polygon and the C shaped polygon, the
RMSE averaged over all the vertices obtained in \cite{elad} when $p$
is known and equal to the true value for
$\sigma=1e^{-3},1e^{-4},1e^{-5}$ is reported in the third column. In
the fourth column the corresponding RMSE obtained by the proposed
procedure is reported.} \label{tb4}
\end{scriptsize}
\end{table}

\subsection{Shape from moments problems}

\begin{figure}
\hspace{.7cm}\epsfig{file=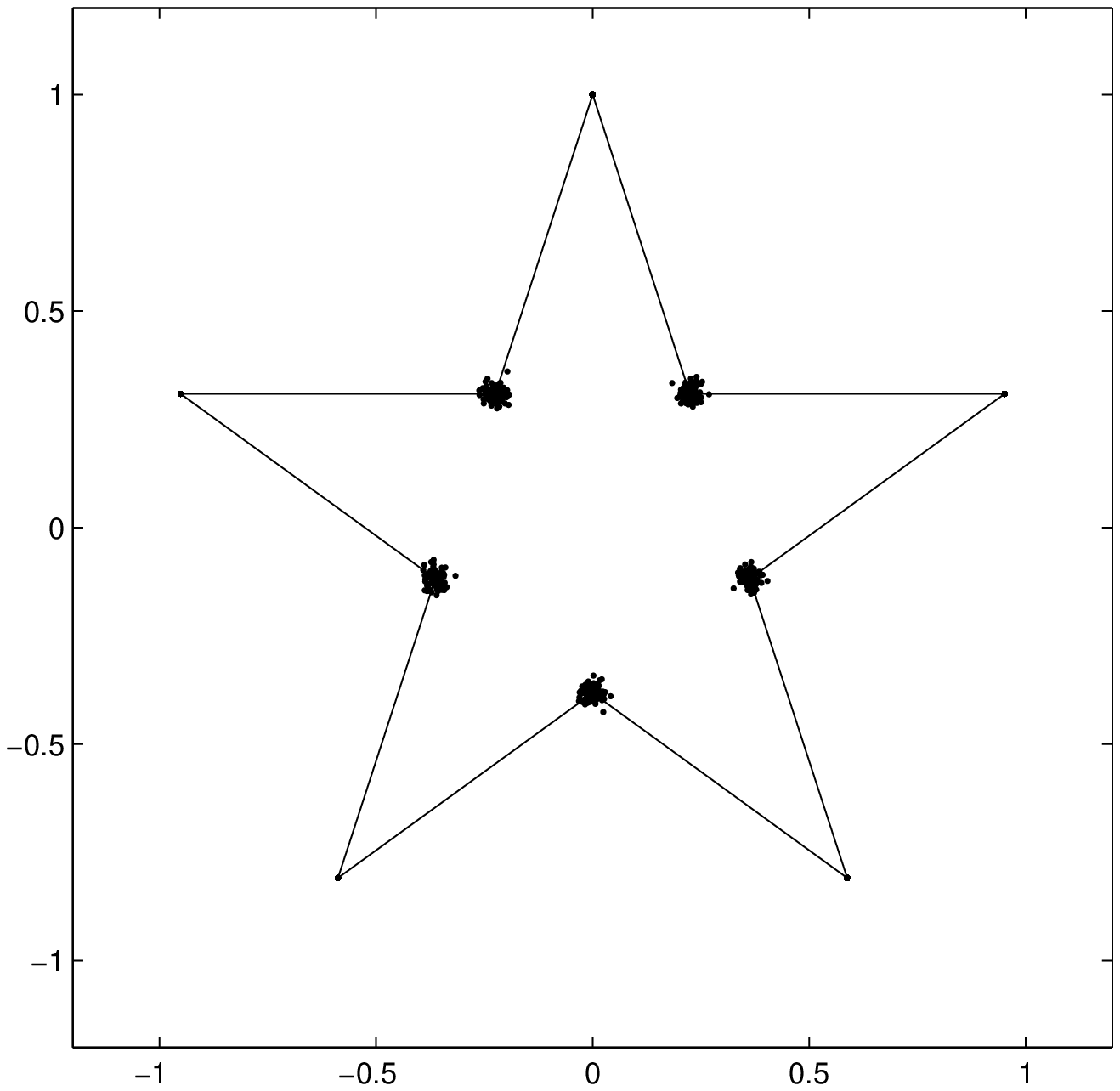,height=6cm}
\hspace{.2cm}\epsfig{file=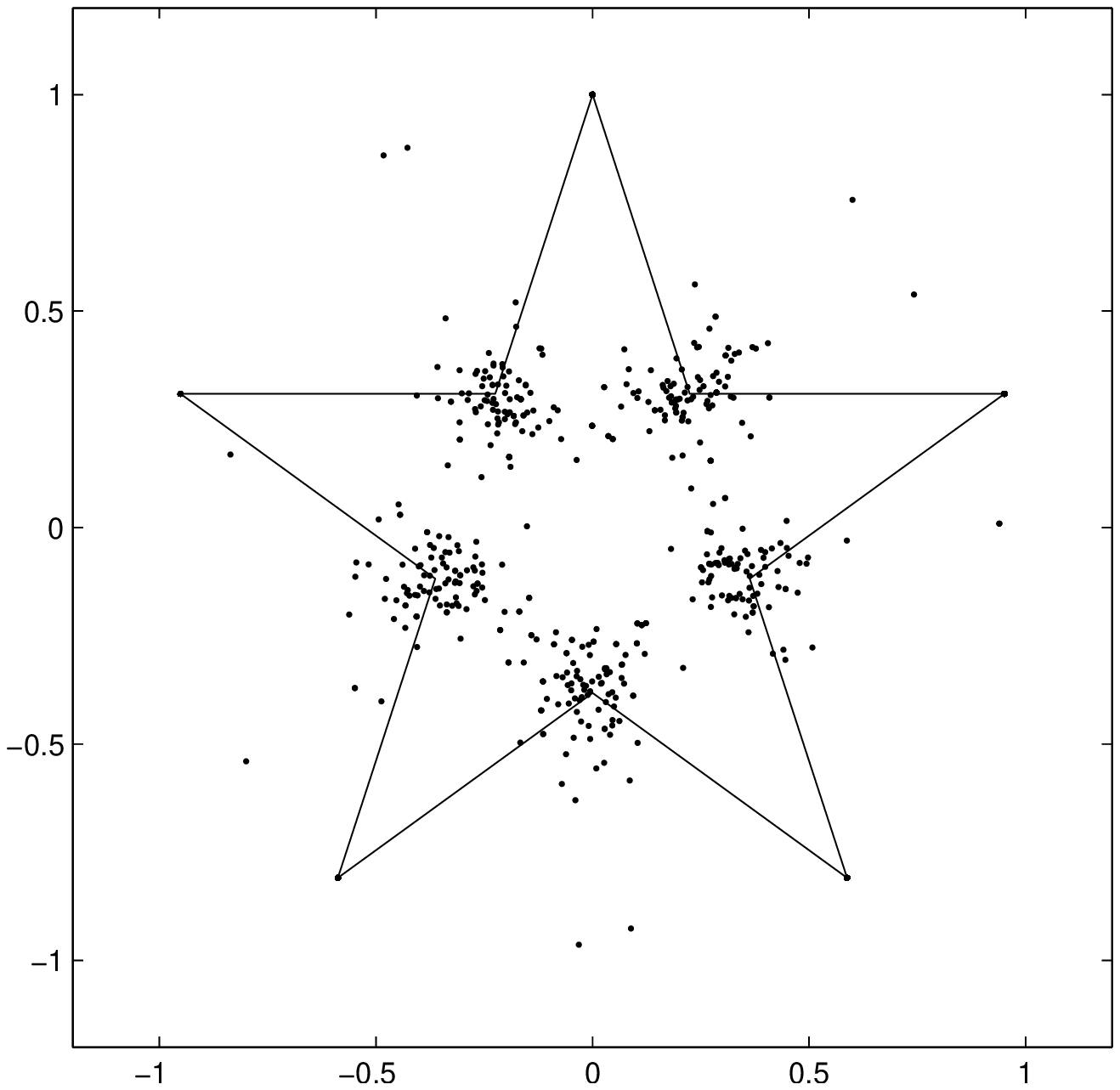,height=6cm}

\vspace{.2cm}

\hspace{.7cm}\epsfig{file=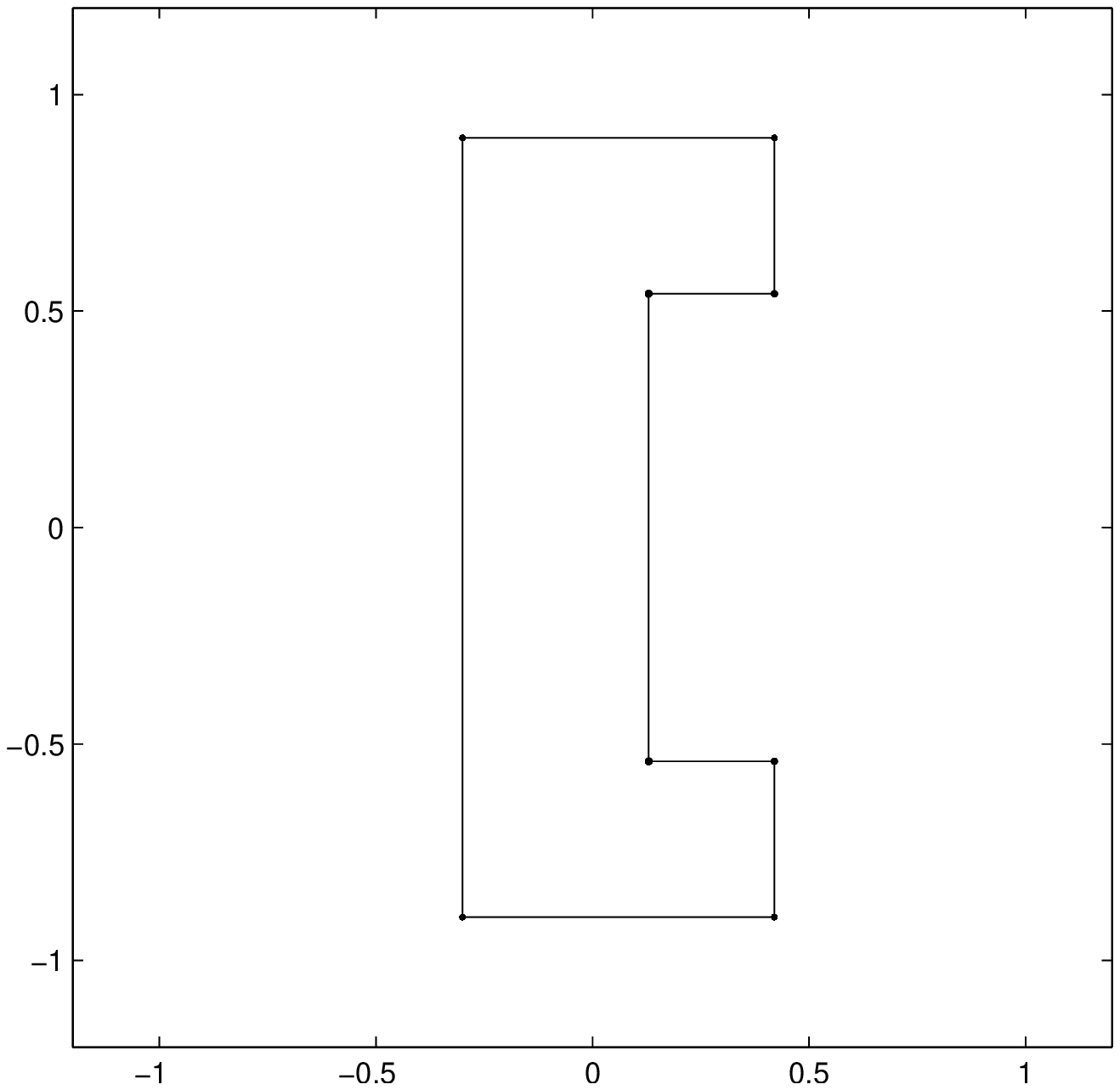,height=6cm}
\hspace{.2cm}\epsfig{file=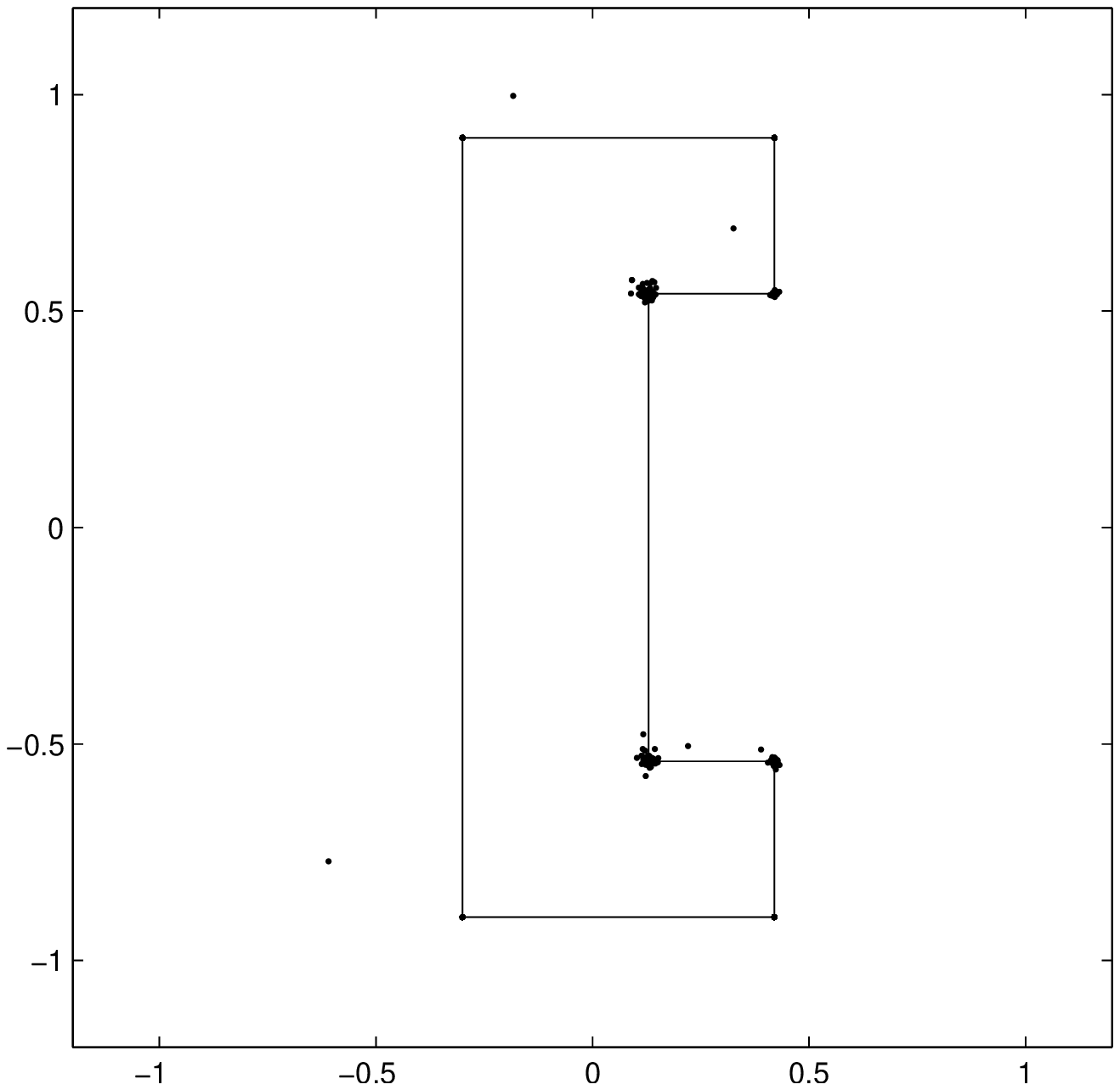,height=6cm}
\caption{Estimates of the vertices of the star shaped and C shaped
polygons obtained by the proposed method on $N=100$ replications
with   $\sigma=1.e^{-4}$. Left: the true number of vertices is
known. Right: the true number of vertices is unknown.}
 \label{fig4}
\end{figure}

In \cite{davis,gmv} it was shown that the $p$ vertices
$\xi_1,\dots,\xi_p$  of a non degenerate polygon $\mathcal{P}$ and
its complex moments $\mu_k,k=0,1,\dots,2p-1$ are related  by
$$k(k-1)\mu_k=k(k-1)\int_{\mathcal{P}}z^k dx\; dy=\sum_{j=1}^p c_j\xi^j,\;\;\mu_0=\mu_1=0$$
where
$$c_j=\frac{i}{2}\left(\frac{\overline{\xi}_{j-1}-\overline{\xi}_j}{\xi_{j-1}-\xi_j}-
\frac{\overline{\xi}_{j}-\overline{\xi}_{j+1}}{\xi_{j}-
\xi_{j+1}}\right)$$ assuming that the vertices are arranged in
counterclockwise direction in the order of increasing index and
extending the indexing of the $\xi_j$ cyclically  so that
$\xi_0=\xi_p$, $\xi_1=\xi_{p+1}$. Therefore to identify the polygon
(i.e. its vertices) from its complex moments is equivalent to solve
a CEIP for the data $s_k=k(k-1)\mu_k$. In \cite{elad} several
methods for solving this specific problem were compared on two
different polygons for $\sigma=10^{-3},10^{-4},10^{-5}$ by a
simulation experiment involving $N=100$ independent replications and
$n=101$ noisy moments. For comparison, in Table 1 the results
obtained by the proposed method and the best among those reported in
\cite[Tables IV, VIII, bold figures]{elad} are reported. The root
mean squared error (RMSE) averaged over all parameters $\xi_j$ is
computed by
$$RMSE=\frac{1}{p}\sum_{k=1}^p\sqrt{ \frac{1}{N}
\sum_{j=1}^{N}|\xi^{(j)}_k-\hat{\xi}^{(j)}_k |^2 }. $$ As the best
results were obtained in \cite{elad} by using GPOF method
(\cite{hua}) but in one case, also in the proposed procedure the
solution of the generalized eigenvalue problem (step 1) was obtained
by GPOF with the same setup used in \cite{elad}. Therefore
$\tilde{p}=p$  is assumed to be known, as in \cite{elad}, and the
$\mathsf{P}$-transform was not computed because all the $p$
estimated clusters were retained. Moreover this was the only example
where GPOF was used also for computing
$\xi_j^{(r)},\;j=1,\dots\tilde{p}, r=1,\dots,R$ (slow method). An
improvement can be noticed in all cases. In the first column of fig.
\ref{fig4} the estimated $\xi_j$ for $\sigma=10^{-4}$ and for the
considered polygons are plotted. We notice that in some vertices,
the $\xi_j$ are so concentrated that they coincide with one point at
the used resolution. Next we use the full fast proposed procedure
assuming not to know $p$ and putting $\tilde{p}=n/2$. The RMSE
averaged over all parameters $\xi_j$ and the mean and standard
deviation of $\hat{p}$ are reported in Table 2. In the second column
of fig. \ref{fig4} the estimated $\xi_j$ for $\sigma=10^{-4}$ and
for the considered polygons are plotted.

\begin{table}[tbh]
\begin{scriptsize}
\begin{center}
\begin{tabular}{|c||c|c|c|c|c|}
\hline\hline
&$\sigma$&$RMSE$& $p$ mean &$p$ s.d. \\
\hline\hline
Star shape&1e-3&1.07e-1&8&3\\
\hline
&1e-4&7.62e-2&10&3\\
\hline
&1e-5&2.98e-2&11&4\\
\hline \hline
C shape&1e-3&4.13e-2&9&1\\
\hline
&1e-4&2.94e-2&9&2\\
\hline
&1e-5&2.72e-2&9&2\\
\hline \hline
\end{tabular}
\end{center}
\caption{For the star shaped polygon and the C shaped polygon, the
RMSE averaged over all the vertices obtained in the case of $p$
unknown for $\sigma=1e^{-3},1e^{-4},1e^{-5}$ is reported in the
third column. In the fourth and fifth columns the estimated mean and
s.d. of $p$ are reported.} \label{tb5}
\end{scriptsize}
\end{table}

\section{Conclusion}

A new approach for solving a classic inverse ill-posed problem is
discussed from the computational point of view. The approach is a
perturbative one, therefore it exploits the information generated by
solving several closed problems by any standard method which best
suits the user's needs such e.g. numerical quality and/or
computational speed. The final results are obtained by an
"averaging" step, hence they are quite stable with respect to noise
and, provided that some hyperparameters are properly selected,
sensitivity is also preserved, allowing to retrieve features of the
signal which are masked by the noise. Several numerical examples are
presented which confirm these practical abilities often improving on
the results given by known methods.

\section*{Acknowledgments}
I wish to thank S.Grande and L.Guidoni  of the Istituto Superiore di
Sanita', Rome, Italy, for providing the NMR data and for many useful
discussions.

\newpage

\begin{thebibliography}{99}
\bibitem{pb06} {\sc Barone, P.} (2008).
  A new transform for solving the noisy complex
exponentials approximation problem, {\em arXiv:0801.1758}.
\bibitem{barja} {\sc Barone, P.} (2005). On the distribution
of poles of Pade' approximants to the Z-transform
 of complex Gaussian white noise, {\em J. Approx.
 Theory} {\bf 132} 224-240.
\bibitem{bama2}  { \sc Barone, P., March, R.} (2001).  A novel class
of Pad\'{e} based method in spectral analysis. {\em J. Comput.
Methods Sci. Eng.} {\bf 1} 185-211.
\bibitem{belk} {\sc Belkic Dz.,  Dando P.A.,  Main J., Taylor H.S. }
(2000). Three novel high-resolution nonlinear methods for fast
signal processing, {\em J.Chem. Phys.} {\bf 113 (16)} 6542--6556
\bibitem{randpol} {\sc Bharucha-Reid A.T., Sambandham M.} (1986).
 {\em Random Polynomials}. Academic Press,
New York.
\bibitem{brez}{\sc Brezinski C., Redivo-Zaglia M.} (1991).
 {\em Extrapolation methods: theory and practice}. North Holland,
Amsterdam.
\bibitem{but}{\sc Boutry, G., Elad, M. Golub, G., Milanfar,
P.} (2005). The generalized eigenvalue problem for nonsquare pencils
using a minimal perturbation approach. {\em SIAM J. Sci. Comp.},{\bf
27,2} 582-601.
\bibitem{cats} {\sc Lendasse, A., Oja, E.,  Simula, O., Verleysen, M.} (2004).  Time Series Prediction Competition:
The CATS Benchmark {\em IJCNN'2004 proceedings – International Joint
Conference on Neural Networks Budapest (Hungary), 25-29 July 2004,
IEEE}, 1615-1620.
\bibitem{davis}{\sc Davis, P.J.} (1964). Triangle formulas in the
complex plane. {\em Math. Comput.}, {\bf 18} 569-577.
\bibitem{donoho}{\sc Donoho, D.L.} (1992).  Superresolution via sparsity constraints.
{\em SIAM J. Math. Anal.}, {\bf 23,5} 1309-1331.
\bibitem{ham}{\sc Hammersley, J.M.} (1956).  The zeros of a random
polynomial. {\em Proc. Berkely Symp. Math. Stat. Probability, 3rd},
{\bf 2} 89-111.
\bibitem{elad}{\sc Elad, M., Milanfar, P., Golub, G.} (2004). Shape from Moments - An
Estimation Theory Perspective, {\em IEEE Trans. on Signal
Processing}, {\bf 52} 1814-1829.
\bibitem{gl} { \sc Golub G.H., Van Loan C.F.} (1996). {\em Matrix computations},
The Johns Hopkins University Press, Baltimore.
\bibitem{gmv}{ \sc Golub, G.H., Milanfar, P., Varah, J.} (2004). A stable numerical method
for inverting shapes from moments. {\em SIAM J. Sci. Comp.},{\bf
21,4} 1222--1243.
\bibitem{hen2}{ \sc Henrici, P.}(1977). { \em Applied and computational
complex analysis, vol.I}, John Wiley, New York.
\bibitem{hua}{\sc Hua, Y., Sarkar, T.K.} (1991). Matrix pencil method
for estimating parameters of damped/undamped sinusoids in noise,
{\em IEEE TASSP},{\bf 39} 892-900.
\bibitem{lem}{\sc Lemmerling, P., Van Huffel, S.} (2002). Structured
total least squares: analysis, algorithms and applications. {\em In
Van Huffel, S., Lemmerling, P.(Eds.) Total least squares and
errors-in-variables modelling}. Kluver, Dordrecht, 79--91.
\bibitem{neu} {\sc Neuhauser, D.} (1990). Bound state eigenfunctions from
wave packets: time-energy resolution, {\em J.Chem. Phys.},{\bf  93}
 2611--2616.
\bibitem{osb}{\sc Osborne M.R.,   Smyth G.K.} (1995).
A Modified Prony Algorithm for Exponential Function Fitting,{ \em
SIAM J. Sci. Comput.} {\bf 16 } 119-138.
 \bibitem{pr} {\sc Prony, R.} (1795). Essai exp\'{e}rimental et analytique sur les lois de
la dilatabilit\'{e} de fluides \'{e}lastiques et sur celles de la
force expansive de la vapeur de l'eau et de la vapeur de l'alkool,
\`{a} diff\'{e}rentes temp\'{e}ratures, { \em Journal de l'\'{E}cole
Polytechnique Flor\'{e}al et Plairial}, III, vol.1, n.22, 24-76.
\bibitem{ps} { \sc Paige, C.C. and Saunders, M.A. }(1982).
  LSQR: An Algorithm for Sparse
Linear Equations And Sparse Least Squares. {\em ACM Trans. Math.
Soft.} {\bf 8} 43-71.
\bibitem{saff}{\sc Saff, E.B., Totik, V.} (1997). {\em Logarithmic
potentials with external fields}, Springer, Berlin
\bibitem{scharf} {\sc Scharf, L.L.} (1991). { \em Statistical signal processing},
Addison-Wesley, Reading.
\bibitem{schu}{ \sc Schuermans M., Lemmerling P., De Lathauwer L., Van Huffel S.} (2006).
 The use of total least squares data fitting in the shape-from-moments problem. {\em
 Signal Processing } {\bf 86} 1109-1115.
\bibitem{fact}{ \sc Sitton, G.A.; Burrus, C.S.; Fox, J.W.;
Treitel, S.} (2003). Factoring very-high-degree polynomials,{\em
Signal Processing Magazine IEEE} {\bf 20}  27- 42.
\bibitem{stw} {\sc Stewart, G.W.} (2001). { \em Matrix algorithms,
 vol.2}, SIAM, Philadelphia.
 \bibitem{huf}{\sc Van Huffel S., Vanderwalle J.} (1991).
 { \em The total least squares problem: computational aspects
 and analysis}, SIAM, Philadelphia.
\bibitem{vin1} {\sc Viti, V., Petrucci, C. and Barone, P.} (1997). Prony methods in NMR
spectroscopy, { \em International Journal of Imaging Systems and
Technology} {\bf 8} 565--571.
\bibitem{vin2} {\sc Viti, V., Ragona, R., Guidoni, G., Barone, P.,
Furman, E., Degani, H.} (1997). Hormonal -induced modulation in the
phosphate metabolites of breast cancer: analysis of in vivo 31P MRS
signals with a modified Prony method, { \em Magnetic Resonance in
Medicine} {\bf 38} 285--295.
 \bibitem{wil}{ \sc Wilkinson J.H.} (1965).{ \em The Algebraic
 Eigenvalue Problem}, Clarendon Press, Oxford.
\end{thebibliography}
\end{document}